\documentclass[letterpaper, 10 pt, conference]{ieeeconf}

\IEEEoverridecommandlockouts

\usepackage{graphicx}
\usepackage{amsmath}
\usepackage{amssymb}
\usepackage[percent]{overpic}
\usepackage[cdot,meduiumqspace,squaren,Gray]{SIunits}
\usepackage{bm}
\usepackage{hyperref}
\usepackage{caption}
\usepackage{subfig}

\newtheorem{theorem}{Theorem}[section]

\DeclareMathOperator*{\ext}{ext}

\title{\LARGE \bf Hamel's Formalism and Variational Integrators on a Sphere}

\author{Dmitry V. Zenkov, Melvin Leok, and Anthony M. Bloch
\thanks{D.\ V.\ Zenkov is with the Department of  Mathematics, North Carolina State University,
      Raleigh, NC 27695, USA,
      {\tt dvzenkov@ncsu.edu}
}%
\thanks{M.\ Leok is with the Department of  Mathematics, University of California San Diego, La Jolla, CA 92093, USA,
{\tt mleok@math.ucsd.edu}
}%
\thanks{A.\ M.\ Bloch is with the Department of  Mathematics,
      University of Michigan, Ann Arbor, MI 48109, USA,
      {\tt abloch@umich.edu}
}%
}

\begin{document}
\maketitle
\thispagestyle{empty}
\pagestyle{empty}
\begin{abstract}
This paper discusses Hamel's formalism and its applications to structure-preserving integration of mechanical systems. It utilizes redundant coordinates in order to eliminate multiple charts on the configuration space as well as nonphysical artificial singularities induced by local coordinates, while keeping the minimal possible degree of redundancy and avoiding integration of differential-algebraic equations.
\end{abstract}

\section{INTRODUCTION}
This paper introduces a new variational integrator for a spherical pendulum. The configuration space for this pendulum is a two-dimensional sphere. Calculations in spherical coordinates are not a good option because of unavoidable artificial singularities introduced by these coordinates at the poles. In addition, the topology of a sphere makes it impossible to use global singularity-free intrinsic coordinates.
 
In order to avoid the issues mentioned above, an integrator that utilizes the interpretation of  a sphere as a homogeneous space was introduced in \cite{LeLeMc2009}. This integrator performs very well, but has a somewhat large degree of redundancy.
This paper targets the development of an integrator whose performance is similar to that of the integrator in \cite{LeLeMc2009} and whose redundancy is the minimum possible.

Both the present paper and \cite{LeLeMc2009} interpret the pendulum as a rotating rigid body. The algorithm introduced in \cite{LeLeMc2009} is based on the evaluation of the rotation matrix that represents the attitude for this body. The key feature of the dynamics utilized in the present paper is that, in order to capture the orientation of a rigid body, it is sufficient to evaluate just one column of that rotation matrix. The resulting equations of motion are interpreted as Hamel's equations written in redundant coordinates.

The general exposition of discrete Hamel's formalism will be a subject of a future publication. Here we demonstrate the usefulness of some of this formalism by constructing an integrator for a spherical pendulum that is energy- and momentum-preserving. The calculations are carried out in the Cartesian coordinates of the three-dimensional Euclidean space. This allows one to avoid singularities and/or multiple coordinate charts that are inevitable for calculations on a sphere. Hamel's approach allows one, among other things, to represent the dynamics in such a way that the length constraint becomes unnecessary. Thus, one avoids the well-known difficulty of numerically solving differential-algebraic equations.

The paper is organized as follows. Hamel's formalism and its discretization are briefly discussed in Sections \ref{mech.sec} and~\ref{disc_mech.sec}. 
The dynamics of a spherical pendulum is reviewed in Section~\ref{pendulum.sec}.
The discrete model for the pendulum based on Hamel's formalism, its comparison to some other discretization techniques, and simulations are given in Sections \ref{discr_pend.sec}, \ref{comparison.sec}, and  \ref{simulations.sec}.


\section{LAGRANGIAN MECHANICS}\label{mech.sec} 
\subsection{The Euler--Lagrange Equations}
A Lagrangian mechanical system is specified by a smooth manifold $Q$ called the \emph{\textbf{configuration space}} and a function $L:TQ\to \mathbb{R} $ called the \emph{\textbf{Lagrangian}}. 
In many cases, the Lagrangian is the kinetic minus potential energy of the system, with the kinetic energy defined by a Riemannian metric and the potential energy being a smooth function on the configuration manifold $Q$.
If necessary, non-conservative forces can be introduced (e.g., gyroscopic forces that are represented by terms in $L$ that are linear in the velocity), but this is not discussed in detail in this paper.

In local coordinates $q=(q^1,\dots,q^n)$ on the configuration space
$Q$ we write $L=L(q,\dot q)$. The dynamics is given by the \textbf{\emph{Euler--Lagrange equations}}
\begin{equation}\label{EL.eqn}
\frac{d}{dt} \frac{\partial L}{\partial \dot{q}^i } = \frac{\partial L}{\partial q^i}, \quad i = 1,\dots,n.
\end{equation} 

These equations were originally derived by Lagrange \cite{La1788} in 1788 by requiring that simple force balance $\bm{F} = m \bm{a}$ be \emph{covariant}, \emph{i.e.}\ expressible in arbitrary
generalized coordinates. A variational derivation of the Euler--Lagrange equations, namely Hamilton's principle (also called the principle of critical action), came later in the work of Hamilton \cite{Ha1834} and \cite{Ha1835} in 1834/35.
For more details, see \cite{Bl2003}, \cite{MaRa1999}, and Theorem~\ref {equivalence.thm} below. 

\subsection{The Hamel Equations}
In this paragraph we briefly discuss the Hamel equations. The exposition follows paper \cite{BlMaZe2009a}.

In many cases the Lagrangian and the equations of motion have a
simpler structure when written using velocity components measured against a frame that is unrelated to system's local configuration coordinates. An example of such a system is the rigid body. 

Let $q = (q^1,\dots, q^n)$ be local coordinates on the configuration space $Q$ and $u_i \in TQ$, $i = 1, \dots,n$, be smooth independent \emph{local} vector fields defined in the same coordinate neighborhood.\footnote{In certain cases, some or all of $u_i$ can be chosen to be \emph{global} vector fields on $Q$.}
The components of $u_i$ relative to the basis $\partial /\partial q^j$ will be denoted $\psi_i ^j$; that is,
\[
u_i(q) = \psi _i^j(q) \frac{\partial}{\partial q^j},
\]
where $i, j = 1,\dots,n$ and where summation on $j$ is understood.

Let $\xi = (\xi^1,\dots, \xi^n) \in \mathbb{R}^n$ be the components of the velocity vector $\dot q \in TQ$ relative to the basis $u_1, \dots, u_n$, \emph{i.e.}, 
\begin{equation}\label{noncommuting.variables.eqn}
\dot q = \xi^i u_i(q);
\end{equation}
then
\begin{equation}\label{l_u.eqn}
l(q, \xi) := L(q, \xi^i u_i (q))
\end{equation}
is the Lagrangian of the system written in the local coordinates $(q, \xi)$ on the tangent bundle $TQ$. The coordinates $(q, \xi)$ are Lagrangian analogues of non-canonical variables in Hamiltonian dynamics. 

Define the quantities $c_{ij}^m (q)$ by the equations
\begin{equation}\label{ccoeff.eqn}
[u_i (q),u_j (q)] = c^m_{ij} (q) u_m (q),
\end{equation}
$i, j, m = 1, \dots, n$. These quantities vanish if and only if the vector fields $u_i (q)$, $i = 1,\dots, n$, commute. Here and elsewhere, $[\cdot,\! \cdot]$ is the Jacobi--Lie bracket of vector fields on $Q$.

Viewing $u_i$ as vector fields on $TQ$ whose fiber components equal $0$ (that is, taking the vertical lift of these vector fields), one defines the directional derivatives $u_i [l]$ for a function $l :TQ\to \mathbb{R}$ by the formula 
\[
u_i [l] = \psi_i^j \frac{\partial l}{\partial q^j}.
\]

The evolution of the variables $(q,\xi)$ is governed by the \textbf{\emph{Hamel equations}} 
\begin{equation}\label{hamel.eqn}
\frac{d}{dt} \frac{\partial l}{\partial\xi^j}
= c^m_{ij}\frac{\partial l}{\partial \xi^m} \xi^i + u_j[l].
\end{equation} 
coupled with equations \eqref{noncommuting.variables.eqn}.
If $u_i = \partial / \partial q^i$, equations \eqref{hamel.eqn} become the Euler--Lagrange equations~\eqref{EL.eqn}. 

Equations \eqref{hamel.eqn} were introduced in 
\cite{Ha1904} (see also \cite{NeFu1972} for details and some history). 

\subsection{Hamilton's Principle for Hamel's Equations}\label{cr_action.sec} 
Let $\gamma: [a,b]\to Q$ be a smooth curve in the configuration space. A \textbf{\emph{variation}} of the curve $\gamma (t)$ is a smooth map $\beta : [a,b] \times [- \varepsilon , \varepsilon] \to Q$ that satisfies the condition $\beta (t, 0) = \gamma (t)$. This variation defines the vector field 
\[
\delta \gamma (t) = \left. \frac{\partial \beta (t,s)}{\partial s} \right|_{s=0} 
\]
along the curve $\gamma (t)$.

\begin{theorem}\label{equivalence.thm} 
\emph{
Let $L:TQ \to \mathbb R$ be a Lagrangian and $l: TQ \to \mathbb R$ be its representation in local coordinates $(q, \xi)$. 
Then, the following statements are equivalent:
\begin{itemize}
\item[{\rm (i)}]
The curve $q(t)$, where $a \leq t \leq b$,  is a critical point of the action functional 
\begin{equation} \label{action.eqn}
\int_a^b L (q, \dot q)\,dt
\end{equation} 
on the space of curves in $Q$ connecting $q_a$ to $q_b$ on the interval $[a,b]$, where we choose variations of the curve $q(t)$ that
satisfy $\delta q(a) = \delta q(b) = 0$.
\item[{\rm (ii)}]
The curve $q(t)$  satisfies the Euler--Lagrange equations~\eqref {EL.eqn}.
\item[{\rm (iii)}]
The curve $(q(t), \xi(t))$ is a critical point of the functional 
\begin{equation} \label{hamel_action.eqn}
\int_a^b l(q,\xi)\,dt
\end{equation} 
with respect to variations $\delta \xi$, induced by the variations $\delta q =  \eta^i u_i (q)$, and given by
\begin{equation}\label{xi_var.eqn}
\delta \xi^k = \dot \eta^k + c_{ij}^k (q)\xi^i \eta^j.
\end{equation}
\item[{\rm (iv)}]
The curve $(q(t), \xi(t))$ satisfies the Hamel equations 
\eqref{hamel.eqn}
coupled with the equations 
$
\dot q = \xi^i u_i (q).
$
\end{itemize}
}
\end{theorem}

\noindent
For the proof of Theorem \ref{equivalence.thm} and the early development and history of these equations, as well as other variational structures associated with Hamel's equations see \cite{Po1901},~\cite{Ha1904},  \cite{BlMaZe2009a}, and \cite{BaZeBl2012}. 

\section{THE DISCRETE HAMEL EQUATIONS}\label{disc_mech.sec}

\subsection{Discrete Hamilton's Principle}\label{disc_hamilton.sec}

A discrete analogue of Lagrangian mechanics can be obtained by discretizing Hamilton's principle; this approach underlies the construction of variational integrators. See Marsden and West~\cite{MaWe2001}, and references therein, for a more detailed discussion of discrete mechanics.

A key notion is that of the \emph{\textbf{discrete Lagrangian}}, which is a map $L^d: Q\times Q \rightarrow \mathbb{R}$ that approximates
the action integral along an exact solution of the Euler--Lagrange
equations joining the configurations $q_k, q_{k+1} \in Q$,
\begin{equation}
\label{exact_ld}
L^d(q_k,q_{k+1})\approx \ext_{q\in\mathcal{C}([0,h],Q)} \int_0^h
L(q,\dot q)\,dt,
\end{equation}
where $\mathcal{C}([0,h],Q)$ is the space of curves
$q:[0,h]\rightarrow Q$ with $q(0)=q_k$, $q(h)=q_{k+1}$, and $\ext$
denotes extremum. 

In the discrete setting, the action integral of
Lagrangian mechanics is replaced by an action sum
\begin{equation*}
      S^d (q_0, q_1, \dots, q_N) = \sum_{k=0}^{N-1}L^d (q_k, q_{k+1}),
\end{equation*}
where $q_k\in Q$, $k = 0, 1, \dots, N$, 
is a finite sequence in the configuration space. 
The equations are obtained by the \emph{\textbf{discrete Hamilton's 
principle}}, which extremizes the discrete action given
fixed endpoints $q_0$ and $q_N$. Taking the extremum over $q_1,\dots
,q_{N-1}$ gives the \emph{discrete Euler--Lagrange equations}
\begin{equation*}
        \label{DEL}
        D_1 L^d (q_k, q_{k+1}) + D_2 L^d (q_{k-1}, q_k) = 0,
\end{equation*}
for $k = 1,\dots ,N-1$. This implicitly defines the update map
$\Phi:Q\times Q\rightarrow Q\times Q$, where
$\Phi(q_{k-1},q_k)=(q_k,q_{k+1})$ and $Q \times Q$ replaces the velocity phase space $TQ$ of Lagrangian mechanics.

In the case that $Q$ is a vector space, it may be convenient to use $(q_{k+1/2}, v_{k,k+1})$, where $q_{k+1/2} = \tfrac12(q _k +q_{k+1})$ and $v_{k,k+1} = \tfrac1h (q_{k+1} - q_k)$, as a state of a discrete mechanical system. In such a representation, the discrete Euler--Lagrange equations become
\begin{align*}
&\tfrac12
\big(
	D_1 L^d (q_{k-1/2}, v_{k-1,k}) + D_1 L^d(q_{k+1/2}, v_{k,k+1})
\big)
\\
&+ \tfrac1h
\big(
	D_2 L^d (q_{k-1/2}, v_{k-1,k}) - D_2 L^d(q_{k+1/2}, v_{k,k+1})
\big) = 0.
\end{align*}
These equations are equivalent to the variational principle
\begin{align}\label{midpoint_principle.eqn}
\nonumber
\delta S^d = \sum_{k=0}^{N-1} &
\big(
D_1 L^d(q_{k+1/2}, v_{k,k+1}) \delta q_{k+1/2}
\\
&+D_2 L^d(q_{k+1/2}, v_{k,k+1}) \delta v_{k,k+1}
\big)
= 0,
\end{align}
where the variations $\delta q_{k+1/2}$ and $\delta v_{k,k+1}$ are induced by the variations $\delta q_k$ and are given by the formulae
\[
\delta q_{k+1/2} = \tfrac12 \big(\delta q_{k+1} + \delta q_k\big),
\quad
\delta v_{k,k+1} = \tfrac1h \big(\delta q_{k+1} - \delta q_k\big).
\]

\subsection{Discrete Hamel's Equations}\label{disc_hamel.sec} 
In order to construct the discrete Hamel equations for a given mechanical system, one starts by selecting the vector fields $u_1(q), \dots, u_n(q)$ and computing the Lagrangian $l(q, \xi)$ given by \eqref{l_u.eqn}. One then discretizes this Lagrangian (we only discuss the mid-point rule here) and obtains
\begin{equation}\label{dlh.eqn}
l^d (q_{k+1/2}, \xi_{k,k+1}) =
h l(q_{k+1/2}, \xi_{k,k+1})
\end{equation}
Note that the discretization in \eqref{dlh.eqn} is carried out \emph{after} writing the continuous-time Lagrangian as a function of $(q, \xi)$.

One of the challenges of discretizing the Hamel equations has been understanding the discrete analogue of the bracket term in \eqref{hamel.eqn}. Until recently, it was only known how to handle this for systems on Lie groups (see e.g.~\cite{BoSu1999a} and \cite{MaPeSh1999}). In the discrete model of a spherical pendulum discussed below, these terms vanish, and we will not discuss the approach to discretize the bracket terms in this paper (details on this topic  can be found in \cite{BaZe2011}).

The analogue of the variational principle \eqref{midpoint_principle.eqn} is obtained by setting
\begin{align*}
\delta q_{k+1/2}^i &= \tfrac12 \psi^i_j (q_{k+1/2}) (\eta^j_{k+1} + \eta^j_k),
\\
\delta \xi_{k,k+1}^i &= \tfrac1h (\eta_{k+1}^i - \eta_k^i) + B_k^i,
\end{align*}
where $B_k$ is the discrete analogue of the bracket term in \eqref{xi_var.eqn}.

The discrete Hamel equations read
\begin{align}\label{disc_hamel.eqn}
\nonumber
&\tfrac12
\big(
	D_u l^d (q_{k-1/2}, \xi_{k-1,k}) + D_u l^d(q_{k+1/2}, \xi_{k,k+1})
\big)
\\
\nonumber
& \hspace{4.5em}
+ B_k^*
+ \tfrac1h
\big(
	D_2 l^d (q_{k-1/2}, \xi_{k-1,k})
\\
& \hspace{9em}
- D_2 l^d(q_{k+1/2}, \xi_{k,k+1})
\big) = 0,
\end{align}
where $D_u l^d$ is the directional derivative given by the formula
\begin{multline*}
D_u l^d(q_{k+1/2}, \xi_{k,k+1})
\\
= \frac{d}{ds} \bigg|_{s=0}  l^d (q_{k+1/2} + s u(q_{k+1/2}),\xi_{k,k+1})
\end{multline*}
and where $B_k^*$ is the discrete analogue of the bracket term of the continuous-time Hamel equations \eqref{hamel.eqn}.
As mentioned above, this term vanishes for the spherical pendulum problem, and is therefore not explicitly shown here (see \cite{BaZe2011} for details).
Equations \eqref{disc_hamel.eqn} along with the discrete analogue of equations \eqref{noncommuting.variables.eqn} define the update map $(q_{k-1/2}, \xi_{k-1,k}) \mapsto (q_{k+1/2}, \xi_{k,k+1})$.

\section{THE SPHERICAL PENDULUM}\label{pendulum.sec}
Consider a spherical pendulum whose length is $r$ and mass is $m$. We view the pendulum as a point mass moving on the sphere of radius $r$ centered at the origin of $\mathbb{R}^3$. The development here is based on the representation
\begin{align}
\label{momentum.eqn}
\dot{\bm M} &= \bm M \times \bm\Omega + \bm T,
\\
\label{reconstruction.eqn}
\dot{\bm\Gamma} &= \bm\Gamma \times \bm\Omega
\end{align}
of the dynamics of a spherical pendulum, where the pendulum is viewed as a rigid body rotating about a fixed point. Here $\bm\Omega$ is the angular velocity of the pendulum, $\bm M$ is its angular momentum, $\bm\Gamma$ is the unit vertical vector (and thus the constraint $\|\bm\Gamma\| = 1$ is imposed), and $\bm T$ is the torque produced by a force acting on the pendulum, all written relative to the body frame.
 Throughout the paper, boldface characters represent three-dimensional vectors.
Note that the projection of\/ $\bm T$ on~$\bm\Gamma$ is zero.
We assume here that the force is conservative, with potential energy $U(\bm\Gamma)$. For the pendulum, $U(\bm\Gamma) = mg \langle \bm a, \bm\Gamma \rangle = mgr \Gamma^3$, where $\bm a$ is the vector from the origin to the pendulum bob.\footnote{All frames in this paper are orthonormal, and thus the dual vectors are interpreted as regular vectors, if necessary.}
Note that the potential energies for forces like gravity are invariant with respect to rotations about~$\bm\Gamma$.

There are two independent components in equation \eqref{momentum.eqn}. We emphasize that this representation, though redundant, eliminates the use of local coordinates on the sphere, such as spherical coordinates.  More details on this appear below.
Spherical coordinates, while being a nice theoretical tool, introduce artificial singularities at the north and south poles. That is, the equations of motion written in spherical coordinates have denominators vanishing at the poles, but this has nothing to do with the physics of the problem and is solely caused by the geometry of the spherical coordinates. Thus, the use of spherical coordinates in calculations is not advisable.

Another important remark is that the length of the vector~$\Gamma$ is a \emph{conservation law} of equations \eqref{momentum.eqn} and \eqref{reconstruction.eqn},
\begin{equation}\label{gamma_length.eqn}
\|\bm\Gamma\| = \mathrm{const},
\end{equation} 
and thus adding the constraint $\|\bm\Gamma\| = 1$ \emph{does not} result in a system of differential-algebraic equations. The latter are known to be a nontrivial object for numerical integration. 

Equations \eqref{momentum.eqn} and \eqref{reconstruction.eqn} may be interpreted in a number of ways. For instance, one can view them as the dynamics of a degenerate rigid body.
For this interpretation, select an orthonormal body frame with the third vector aligned along the direction of the pendulum. The inertia tensor relative to such a frame is\/ $\mathcal{I} = \operatorname{diag}\{mr^2, mr^2, 0\}$, and 
the Lagrangian reads
\begin{equation}\label{pendulum_lagr.eqn} 
l(\bm\Omega, \bm\Gamma) = \tfrac12 \langle \operatorname{\mathcal{I}} \bm\Omega, \bm\Omega \rangle
- U(\bm\Gamma).
\end{equation}
With this frame selection, the third component of the angular momentum of the body vanishes,
\[
M_3 = \frac{\partial l}{\partial \Omega^3} = \mathcal{I}_3 \Omega^3 = 0,
\]
and thus there are only two nontrivial equations in \eqref{momentum.eqn}. 
Thus, one needs five equations to capture the pendulum dynamics. This reflects the fact that rotations about the direction of the pendulum have no influence on the pendulum's motion. The dynamics then can be simplified by setting $\Omega^3 = 0$.

Alternatively, \eqref{momentum.eqn} and \eqref{reconstruction.eqn} may be interpreted as the dynamics of the \emph{Suslov problem} (see \cite{NeFu1972} and \cite{Bl2003}) for a rigid body with a rotationally-invariant inertia tensor and constraint $\Omega ^3 = 0$.

Using either interpretation, the dynamics is represented by the system of five first-order differential equations
\begin{align}
\label{ct_momentum.eqn} 
\dot M_1 &= mgr \Gamma^2,
\quad
\dot M_2 = -mgr \Gamma^1,
\\
\label{ct_reconstruction.eqn}
\dot \Gamma^1 &= - \Omega^2 \Gamma^3,
\quad
\dot \Gamma^2 = \Omega^1 \Gamma^3,
\quad
\dot \Gamma^3 = \Omega^2 \Gamma^1 - \Omega^1 \Gamma^2.
\end{align}

The latter equations are in fact Hamel's equations written in the redundant coordinates $(\Gamma^1, \Gamma^2, \Gamma^3)$ relative to the frame
\begin{align*}
u_1 &=
\Gamma^3 \frac{\partial }{\partial \Gamma^2} - 
\Gamma^2 \frac{\partial }{\partial \Gamma^3},
\\
u_2 &=
\Gamma^1 \frac{\partial }{\partial \Gamma^3} - 
\Gamma^3 \frac{\partial }{\partial \Gamma^1}.
\end{align*}
Recall that the length of $\bm\Gamma$ is the conservation law, so that the constraint $\|\bm\Gamma\| = 1$ does not need to be imposed, but the appropriate level set of the conservation law needs to be selected.

Our discretization will be based on this point of view, i.e., the discrete dynamics will be written in the form of discrete Hamel's equations.
The discrete dynamics will posses the discrete version of the conservation law \eqref{gamma_length.eqn}, so that the algorithm should be capable, in theory, of preserving the length of $\bm\Gamma$ up to machine precision.

\section{VARIATIONAL DISCRETIZATION FOR THE SPHERICAL PENDULUM}\label{discr_pend.sec}
The integrator for a spherical pendulum is constructed by discretizing Hamel's equations \eqref{ct_momentum.eqn} and \eqref{ct_reconstruction.eqn}.

Let the positive real constant $h$ be the time step. Applying the mid-point rule to \eqref{pendulum_lagr.eqn}, the discrete Lagrangian is computed to be
\begin{equation}\label{disc_lagr_pend.eqn}
l^d = \tfrac{h}{2}
\langle
	\operatorname{\mathcal I} \bm\Omega_{k,k+1}, \bm\Omega_{k,k+1}
\rangle
- h U(\bm\Gamma _{k+1/2}).
\end{equation}
Here $\bm\Omega_{k,k+1} = (\bm\Omega^1_{k,k+1}, \bm\Omega^2_{k,k+1}, 0)$ is the discrete analogue of the angular velocity $\bm\Omega = (\Omega^1, \Omega^2, 0)$ and $\bm\Gamma_{k+1/2} = \frac12 (\bm\Gamma_{k+1} + \bm\Gamma_k)$.
The discrete dynamics then reads
\begin{align}
\label{disc_mom.eqn}
\tfrac1h \operatorname{\mathcal I} \big(\bm\Omega_{k,k+1} - \bm\Omega_{k-1,k}\big)
&= \bm T_{k+1/2},
\\
\nonumber
\label{disc_rec.eqn}
\tfrac1h \big(\bm\Gamma_{k+1/2}  - \bm\Gamma_{k-1/2}\big)
&= \tfrac12 \big(\bm\Gamma_{k+1/2} + \bm\Gamma_{k-1/2}\big)
\\
& \quad\ 
\times 
\tfrac12 \big(\bm\Omega_{k,k+1} + \bm\Omega_{k-1,k}\big),
\end{align}
or, in components,
\begin{align}
\label{disc_mom_1.eqn} 
\tfrac1h r^2 \big(\bm\Omega_{k,k+1}^1 - \bm\Omega_{k-1,k}^1\big) &=
\tfrac12 gr\big(\bm\Gamma_{k+1/2}^2 + \bm\Gamma_{k-1/2}^2\big),
\\
\label{disc_mom_2.eqn} 
\tfrac1h r^2 \big(\Omega_{k,k+1}^2 - \Omega_{k-1,k}^2\big) &= -
\tfrac12 gr\big(\Gamma_{k+1/2}^1 + \Gamma_{k-1/2}^1\big),
\\
\nonumber
\tfrac1h \big(\Gamma_{k+1/2}^1 - \Gamma_{k-1/2}^1\big) &= -
\tfrac14 \big(\Omega_{k,k+1}^2 + \Omega_{k-1,k}^2\big)
\\
\nonumber
& \qquad\, \times 
\big(\Gamma_{k+1/2}^3 + \Gamma_{k-1/2}^3\big),
\\
\nonumber
\tfrac1h \big(\Gamma_{k+1/2}^2 - \Gamma_{k-1/2}^2\big) &=
\tfrac14 \big(\Omega_{k,k+1}^1 + \Omega_{k-1,k}^1\big)
\\
\nonumber
& \qquad\, \times 
\big(\Gamma_{k+1/2}^3 + \Gamma_{k-1/2}^3\big),
\\
\label{disc_rec_3.eqn}
\nonumber
\tfrac1h \big(\Gamma_{k+1/2}^3 - \Gamma_{k-1/2}^3\big) &=
\tfrac14 \big(\Omega_{k,k+1}^2 + \Omega_{k-1,k}^2\big)
\\
\nonumber
& \quad\ \ \times 
\big(\Gamma_{k+1/2}^1 + \Gamma_{k-1/2}^1\big)
\\
\nonumber
&\quad\ \ \ -
\tfrac14 \big(\Omega_{k,k+1}^1 + \Omega_{k-1,k}^1\big)
\\
& \qquad\ \ \times 
\big(\Gamma_{k+1/2}^2 + \Gamma_{k-1/2}^2\big).
\end{align}

The discretization of equation \eqref{reconstruction.eqn} is constructed to be $\|\bm\Gamma\|$-preserving. Indeed, using the isomorphism $\bm\Omega \mapsto \Omega$ between the spaces of three-dimensional vectors $\bm\Omega = (\Omega^1, \Omega^2, \Omega^3)$ and skew-symmetric matrices
\[
\Omega =
\begin{bmatrix}
\hphantom{-}0\hphantom{^3} & - \Omega^3 &\hphantom{-}\Omega^2
\\
\hphantom{-}\Omega^3 & \hphantom{-}0\hphantom{^3} & -\Omega^1
\\
-\Omega^2 & \hphantom{-}\Omega^1 &\hphantom{-} 0\hphantom{^1}
\end{bmatrix}\!,
\]
equation \eqref{disc_rec.eqn} becomes
\[
\bm\Gamma_{k+1/2} =
(I - A_k)^{-1} (I + A_k) \bm\Gamma_{k-1/2},
\]
where
\[
A_k = \tfrac h4 \big(\Omega_{k,k+1} + \Omega_{k-1,k}\big).\footnote{The matrix $I - A_k$ is invertible if $h$ is sufficiently small.}
\]
It is straightforward to check that the matrix
\[
(I - A_k)^{-1} (I + A_k)
\]
is orthogonal (it is simply the Cayley transform of $A_k$), and therefore $\|\bm\Gamma_{k+1}\| = \|\bm\Gamma_k\|$. 

As expected, the discrete dynamics is momentum-preserving: 
\[
m r^2
\big(
\Gamma_{k+1/2}^1 \Omega_{k,k+1}^1
+ \Gamma_{k+1/2}^2 \Omega_{k,k+1}^2
\big)
= \mathrm{const},
\]
i.e., the vertical component of spatial momentum is conserved. This can be verified either using  general symmetry arguments, or by a straightforward calculation.

The discrete dynamics is energy-preserving:
\[
\tfrac12 m r^2
\Big[
\big(\Omega_{k,k+1}^1\big)^2 + \big(\Omega_{k,k+1}^2\big)^2
\Big]
+ mgr \Gamma_{k+1/2}^3 = \mathrm{const}.
\]
This is confirmed by multiplying equations \eqref{disc_mom_1.eqn} and \eqref{disc_mom_2.eqn} by
$\big(\Omega_{k,k+1}^1 + \Omega_{k-1,k}^1\big)$
and 
$\big(\Omega_{k,k+1}^2 + \Omega_{k-1,k}^2\big)$,
respectively, and adding the result to equation \eqref{disc_rec_3.eqn}.

To recap, the proposed method preserves the symplectic structure, the length constraint, and the momentum, so by a result due to Ge and Marsden \cite{GeMa1988}, the proposed method recovers the exact trajectory, up to a possible time reparameterization.

\section{COMPARISON WITH OTHER METHODS}\label{comparison.sec}
The proposed method takes advantage of the homogeneous space structure of $S^2$, which has a transitive Lie group action by $\operatorname{SO(3)}$. In particular, the vector $\Gamma\in S^2$ is updated by the left action of a rotation matrix, given by the Cayley transformation of a skew-symmetric matrix $A_k$ that approximates the angular momentum $\Omega$ integrated over a half-timestep. Interestingly, this falls out naturally from discretizing the Hamel formulation of the spherical pendulum, and it would be interesting to see what general choices of coordinate frames in the Hamel formulation lead to similar methods for more general flows on homogeneous spaces. We will now discuss some alternative methods of simulating the spherical pendulum equations.

\subsection*{Homogeneous Space Variational Integrators}
If one were to instead formulate the spherical pendulum problem directly on $S^2$, it is possible to lift the variational principle on $S^2$ to $\operatorname{SO(3)}$, by relating the curve $\bm\Gamma(t)\in S^2$ with a curve $R(t)\in \operatorname{SO(3)}$, by the relation $\bm\Gamma(t)=R(t)\bm\Gamma(0)$, where $R(0)=I$, except that the resulting variational principle on $\operatorname{SO(3)}$ does not have a unique extremizer, due to the presence of a nontrivial isotropy subgroup. With a suitable choice of connection, this ambiguity can be eliminated, and the resulting problem (and similar problems on homogeneous spaces) can be solved using Lie group variational integrator techniques, as described in \cite{LeLeMc2009}.

\subsection*{Nonholonomic Integrators}
As mentioned in Section \ref{pendulum.sec}, the spherical pendulum equations can be viewed as a Suslov problem, which is an example of a nonholonomic mechanical system with no shape space. In principle, one could apply a nonholonomic integrator, such as the one described in \cite{FeZe2005} and \cite{McPe2006}, but replacing the length constraint with an infinitesimal constraint and a discrete nonholonomic constraint may result in poor numerical preservation of the constraint properties if the discrete nonholonomic constraint is poorly chosen. An alternative approach to simulating nonholonomic mechanics involves a discretization of the forces of constraint, and a careful choice of force discretization has been shown to yield promising results, see \cite{LyZe2009} for details.

\subsection*{Constrained Symplectic Integrators}
Given the relatively simple nature of the unit length constraint, it is quite natural to apply the RATTLE algorithm~\cite{An1983}, which is a generalization of the St\"ormer--Verlet method for constrained Hamiltonian systems that is designed to explicitly preserve holonomic constraints. This method does require the use of a nonlinear solver on a system of nonlinear equations of dimension equal to the number of constraints. The cost of the nonlinear solve can increase significantly as the number of copies of the sphere in the configuration space increases.

\subsection*{Differential-Algebraic Equation Solvers}
The proposed discrete Hamel integrator can be easily scaled to an arbitrary number of copies of the sphere, possibly chained together in a $n$-spherical pendulum. Such multi-body systems however pose significant challenges for differential-algebraic equation solvers, since they are examples of what are referred to as high-index DAEs, for which the theory and numerical methods are much less developed. It is possible to perform index reduction on the system of differential-algebraic equations, but this involves significant effort, and the numerical results can be mixed.

\subsection*{Numerical Comparisons}
Since the method is a second-order accurate symplectic method, it is natural to compare it to the St\"ormer--Verlet method, as well as the RATTLE method (which is a generalization of St\"ormer--Verlet for constrained Hamiltonian systems).

For the St\"ormer--Verlet method, we compute Hamilton's equations for the spherical pendulum, which is given by
\begin{align}
\dot{\bm x} &= \tfrac1m \bm p,\label{index_red1.eq}
\\
\dot{\bm p} &= -mg\bm e_3 + \left(mg\bm x\cdot\bm e_3
-\tfrac1m \|\bm p\|^2\right)\bm x = \bm f(\bm x,\bm p),\label{index_red2.eq}
\end{align}
and we apply the generalization of the St\"ormer--Verlet method for general partitioned problems (see (3.4) in \cite{HaLuWa2006}),
\begin{align*}
\bm p_{n+1/2} &= \bm p_n + \tfrac{h}{2}\bm f(\bm x_n, \bm p_{n+1/2}),
\\
\bm x_{n+1} &= \bm x_n + \tfrac hm \bm p_{n+1/2},
\\
\bm p_{n+1} &= \bm p_{n+1/2} + \tfrac h2 \bm f (\bm x_{n+1}, \bm p_{n+1/2}).
\end{align*}
This system of equations is linearly implicit, since the first equation is implicit in $\mathbf{p}_{n+\frac{1}{2}}$, but the rest of the equations are explicit.

The RATTLE method (see (1.26) in \cite{HaLuWa2006}) can be applied to the particle in a uniform gravitational field problem,
\begin{align*}
\dot{\bm x} &= \tfrac1m \bm p,
\\
\dot{\bm p} &= -mg\bm e_3,
\end{align*}
subject to the constraint $\phi(\bm x) = \frac12 (\|\bm x\|^2-1)=0$. We also introduce $\bm\Phi(\bm x) = \frac{\partial \phi}{\partial \bm x} = \bm x^T$. Then, the RATTLE method applied to this problem is given by,
\begin{align*}
\bm p_{n+1/2} &= \bm p_n - \tfrac h2 \big(mg\bm e_3 + \bm\Phi(\bm x_n)^T\lambda_n\big),
\\
\bm x_{n+1} &= \bm x_n + \tfrac hm \bm p_{n+1/2},
\\
0 &= \phi(\bm x_{n+1}),
\\
\bm p_{n+1} &= \bm p_{n+1/2} - \tfrac h2 \big(mg\bm e_3 + \bm\Phi(\bm x_{n+1})^T\mu_n\big),
\\
0& = \tfrac1m \bm\Phi(\bm x_{n+1})\cdot \bm p_{n+1}.
\end{align*}

\section{SIMULATIONS}\label{simulations.sec}

\begin{figure}[h!]
\begin{center}
\subfloat[Trajectory of the pendulum on $S^2$]{\label{trajectory.fig}\includegraphics[width=.8\columnwidth]{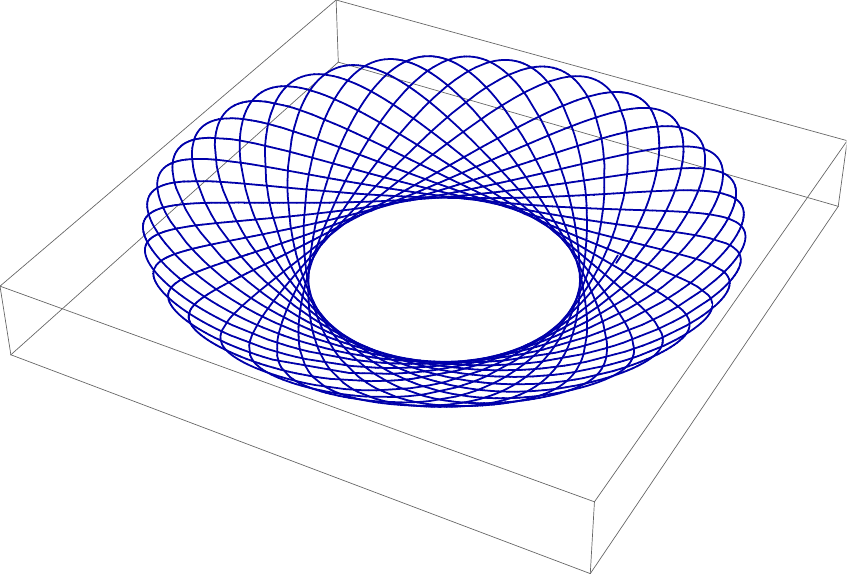}}\\
\subfloat[Preservation of the length of $\bm\Gamma$]{\label{norm.fig}\includegraphics[width=.9\columnwidth]{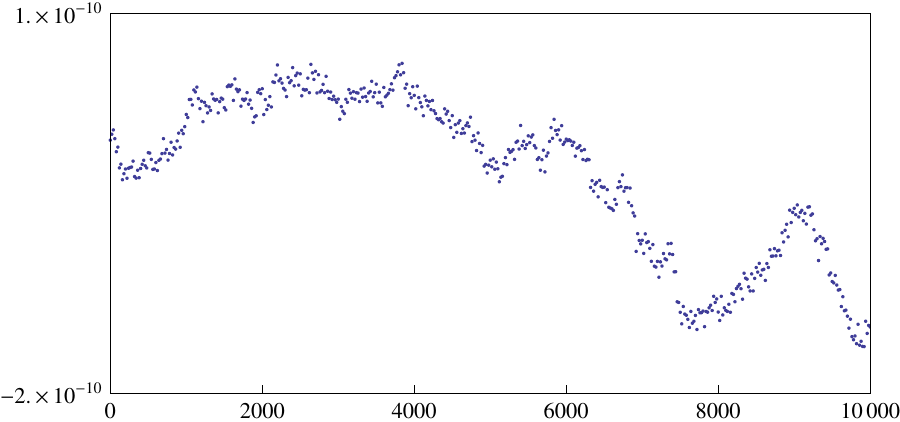}}\\
\subfloat[Conservation of energy]{\label{energy.fig}\includegraphics[width=.9\columnwidth]{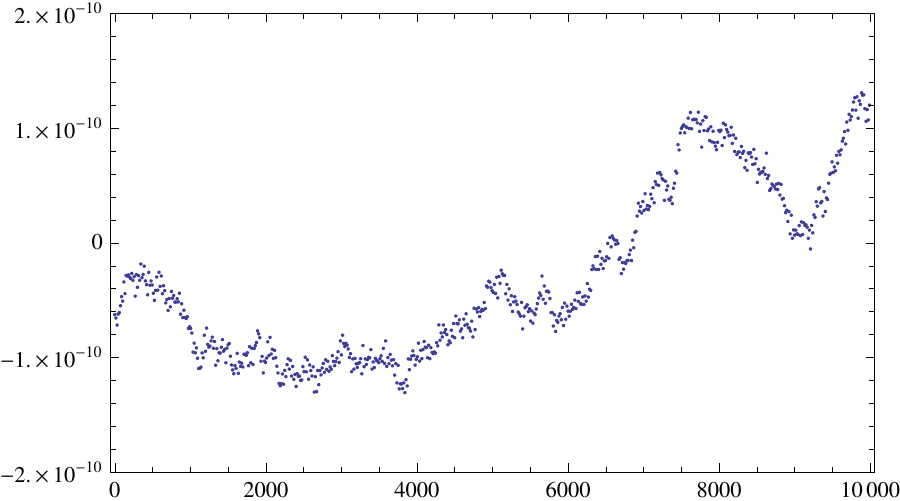}}
\caption{Numerical properties of the Hamel integrator.}
\label{HVI.fig}
\end{center}
\end{figure}

\begin{figure}[h!]
\begin{center}
\includegraphics[width=.45\columnwidth]{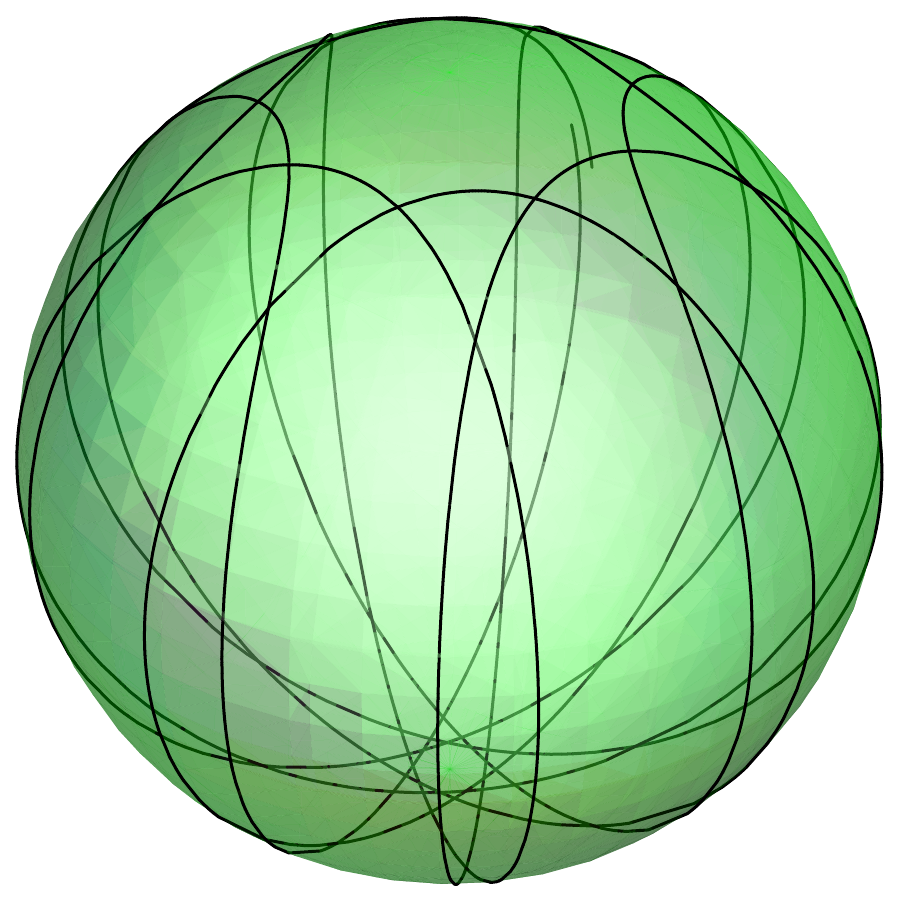}
\caption{A trajectory with initial conditions above the equator integrated with the Hamel integrator.}\label{trajectory2.fig} 
\end{center}
\end{figure}

\begin{figure}[h!]
\begin{center}
\subfloat[Trajectory of the pendulum on $S^2$]{\label{SV_trajectory.fig}\includegraphics[width=.8\columnwidth]{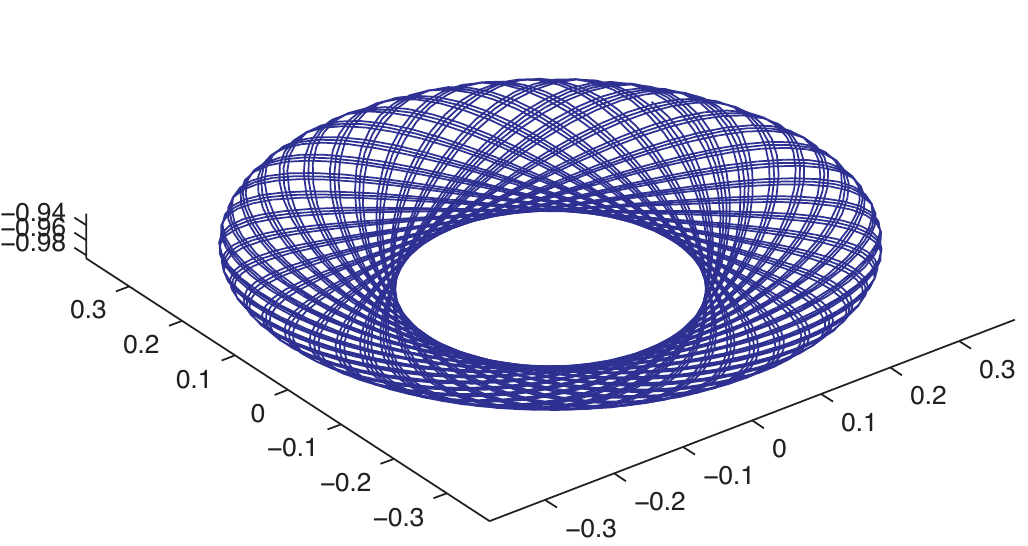}}\\
\subfloat[Preservation of the length of $\bm\Gamma$]{\label{SV_norm.fig}\includegraphics[width=.9\columnwidth]{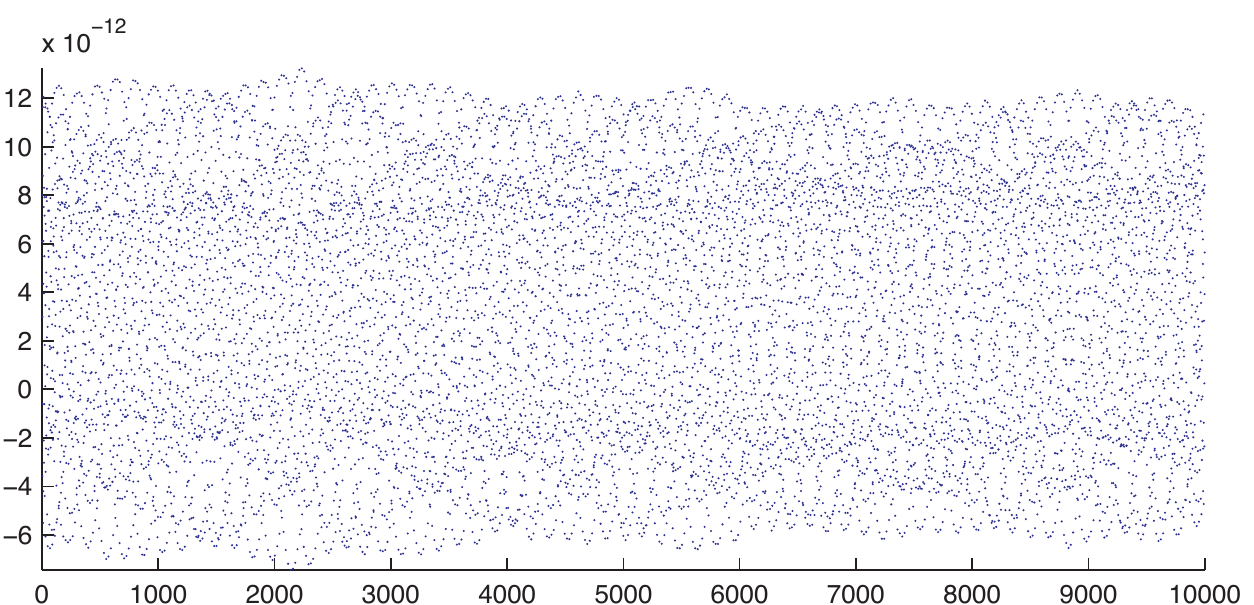}}\\
\subfloat[Conservation of energy]{\label{SV_energy.fig}\includegraphics[width=.9\columnwidth]{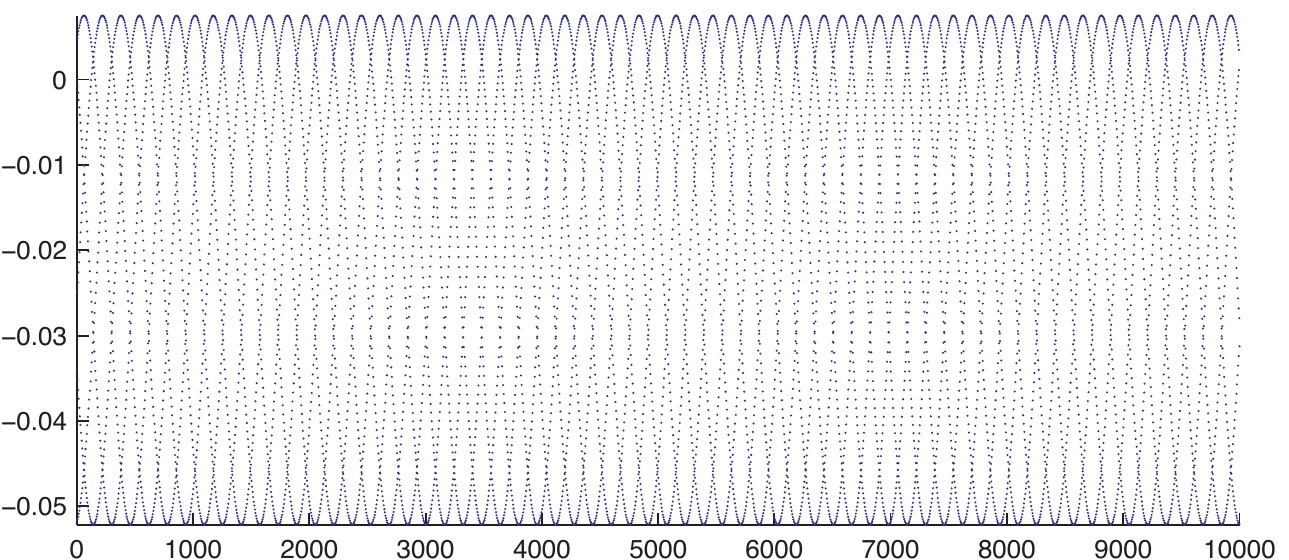}}
\caption{Numerical properties of the St\"ormer--Verlet method.}
\label{SV.fig}
\end{center}
\end{figure}

\begin{figure}[h!]
\begin{center}
\subfloat[Trajectory of the pendulum on $S^2$]{\label{RATTLE_trajectory.fig}\includegraphics[width=.8\columnwidth]{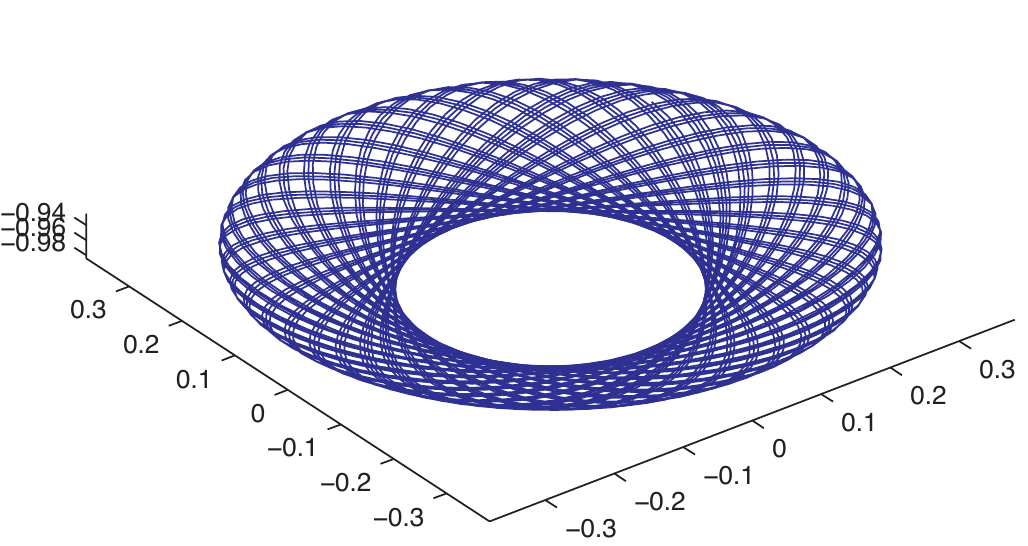}}\\
\subfloat[Preservation of the length of $\bm\Gamma$]{\label{RATTLE_norm.fig}\includegraphics[width=.9\columnwidth]{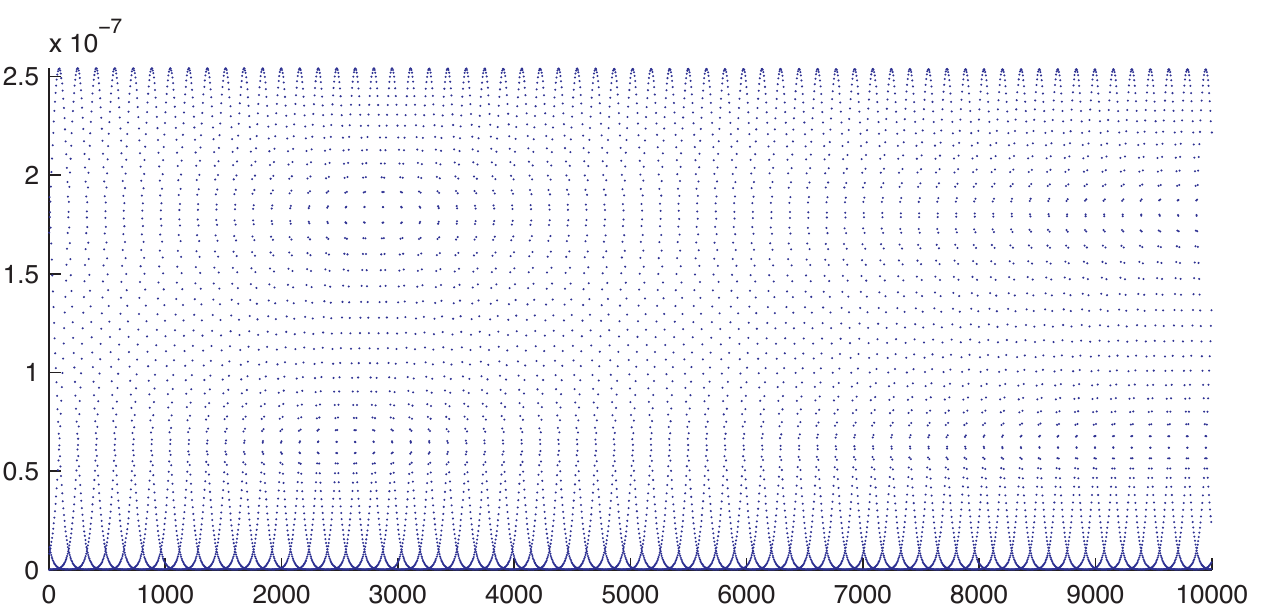}}\\
\subfloat[Conservation of energy]{\label{RATTLE_energy.fig}\includegraphics[width=.9\columnwidth]{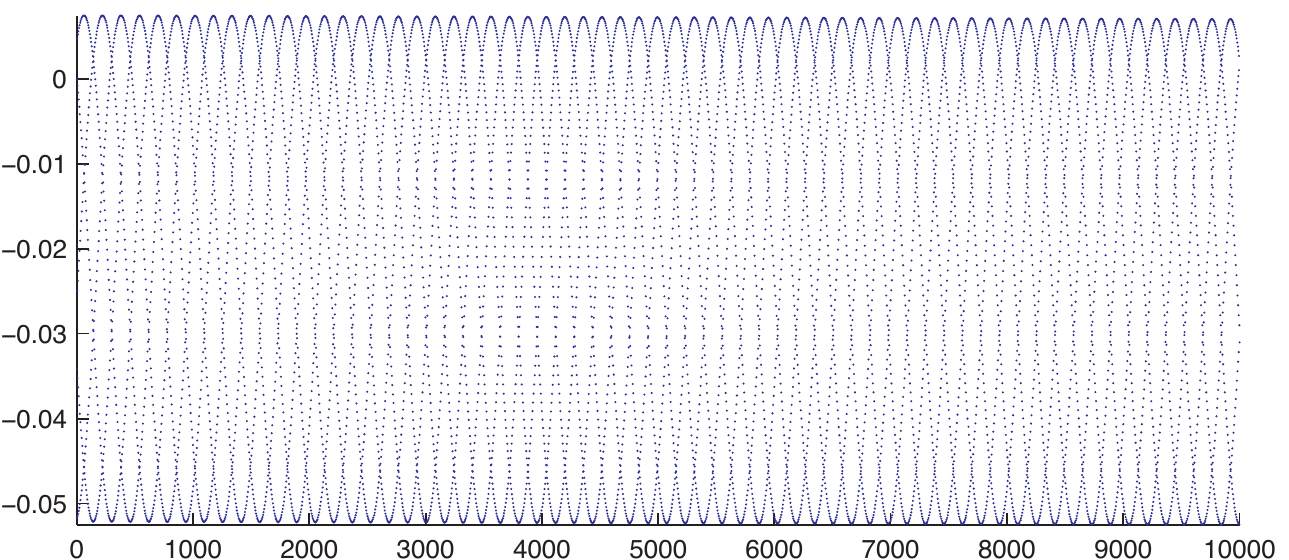}}
\caption{Numerical properties of the RATTLE method.}
\label{RATTLE.fig}
\end{center}
\end{figure}

In Figures \ref{HVI.fig} and \ref{trajectory2.fig}, we present simulations using our theory developed above, which we compare with simulations using the generalized St\"ormer--Verlet method and the RATTLE method in Figures \ref{SV.fig} and \ref{RATTLE.fig}, respectively.

For simulations, we select the parameters of the system and the initial conditions to be 
$m = \unit{1}{\kilogram}$,
$r = \unit{9.8}{\metre}$, $h= \unit{.2}{\second}$,
$\Omega^1_0 = \unit{.6}{\radianpersecond}$,
$\Omega^2_0 = \unit{0}{\radianpersecond}$,
$\Gamma^1_0 = \unit{.3}{\metre}$,
$\Gamma^2_0 = \unit{.2}{\metre}$,
$\Gamma^3_0 = \unit{-.932738}{\metre}$.
The trajectory of the bob of the pendulum is shown in Figure \ref{trajectory.fig}. As expected, it reveals the quasiperiodic nature of pendulum's dynamics.
Theoretically, if one solves the nonlinear equations exactly, and in the absence of numerical roundoff error, the Hamel variational integrator should exactly preserve the length constraint, and the energy. In practice, Figure \ref{norm.fig} demonstrates
that $\|\bm\Gamma\|$ stays to within unit length to about $10^{-10}$ after 10,000 iterations. Figure \ref{energy.fig} demonstrates numerical energy conservation, and the energy error is to about $10^{-10}$ after 10,000 iterations as well. Indeed, one notices that the energy error tracks the length error of the simulation, which is presumably due to the relationship between the length of the pendulum and the potential energy of the pendulum. The drift in both appear to be due to accumulation of numerical roundoff error, and could possibly be reduced through the use of compensated summation techniques.

Figure \ref{trajectory2.fig} shows pendulum's trajectory that crosses the equator.  
This simulation demonstrates the global nature of the algorithm, and also seems to do a good job of hinting at the geometric conservation properties of the method.

\addtolength{\textheight}{-.29cm}

We also simulate the spherical pendulum using the generalized St\"ormer--Verlet method and the RATTLE method described in Section \ref{comparison.sec}. The generalized St\"ormer--Verlet method exhibits surprisingly good unit length preservation in Figure \ref{SV_norm.fig} of $10^{-11}$ when applied to index-reduced version of the equations of motion \eqref{index_red1.eq}--\eqref{index_red2.eq}.
The energy behavior in Figure \ref{SV_energy.fig} is typical of a symplectic integrator, with the characteristic bounded energy oscillations. Even though the RATTLE method is intended to explicitly enforce the unit length constraint, it exhibits a unit length preservation in Figure \ref{RATTLE_norm.fig} of $10^{-7}$, which is poorer than both the Hamel variational integrator and the generalized St\"ormer--Verlet method. The energy error for RATTLE in Figure \ref{RATTLE_energy.fig} is comparable to that of the generalized St\"ormer--Verlet method, but both pale in comparison to the energy error for the Hamel variational integrator.

\vspace{-.35em}
\section{CONCLUSIONS AND FUTURE WORK}
In this paper, we constructed a variational integrator for the spherical
pendulum by discretizing Hamel's equations. We showed the integrator
preserves key mechanical quantities and illustrated the work with simulations, and comparisons with the generalized St\"ormer--Verlet method and the RATTLE method.

Future work will include the simulation of linked rigid body systems as
well as the control of such systems. The excellent numerical properties of the proposed Hamel variational integrator will serve as an excellent basis for constructing numerical optimal control algorithms, which are heavily dependent on the quality of the numerical discretization of the natural dynamics.


\section{ACKNOWLEDGMENTS}
The authors would like to thank the reviewers for valuable comments.
The research of AMB was partially supported by NSF grants DMS-0806765,  DMS-0907949 and DMS-1207993.
The research of ML was partially supported by NSF grants DMS-1010687, CMMI-1029445, and DMS-1065972.
The research of DVZ was partially supported by NSF grants DMS-0306017, DMS-0604108, DMS-0908995, and DMS-1211454. 


\end{document}